\documentclass{article}
\usepackage{epsfig}
\usepackage{amsmath, amssymb, amsthm}	
\usepackage{color}			
\usepackage{graphicx}
\usepackage{caption}
\usepackage{subcaption}

\renewcommand{\vec}[1]{\mbox{\boldmath $ #1 $}}

\newtheorem{remark}{Remark}[section]
\newtheorem{example}{Example}[section]

\begin{document}

\pagestyle{headings}

\title{Comparison of Integrators for Electromagnetic Particle in Cell Methods: Algorithms and Applications}
\author{J\"urgen Geiser\thanks{Ruhr University of Bochum, Department of Electrical Engineering and Information Technology, Universit\"ats­str. 150, D-44801 Bochum, Germany, E-mail: juergen.geiser@ruhr-uni-bochum.de}
  \and Frederik Riedel\thanks{Ruhr University of Bochum, Department of Physics, Universit\"ats­str. 150, D-44780 Bochum, Germany, E-mail: frederik.riedel@ruhr-uni-bochum.de}
}
\maketitle

\begin{abstract}

In this paper, we present different types of integrators
for electro-magnetic particle-in-cell (PIC) methods.
While the integrator is an important tool of the 
PIC methods, it is necessary to characterize the different
conservation approaches of the integrators, e.g. symplecticity,
energy- or charge-conservation.
We discuss the different principles, e.g. composition,
filtering, explicit and implicit ideas.
 
While, particle in cell methods are well-studied,
the combination between the different parts, i.e. pusher,
solver and approximations are hardly to analyze.
we concentrate on choosing the optimal pusher component,
with respect to conservation and convergence behavior.

We discuss oscillations of the pusher component, strong external magnetic
fields and optimal conservation of energy and momentum.

The algorithmic ideas are discussed and numerical experiments
compare the exactness of the different schemes.

An outlook to overcome the different error components is discussed in the
future works.

\end{abstract}



\section{Introduction}

We motivate our studies on simulations of an electro-magnetic Particle-in-Cell
(EMPIC) method and an optimal combination of the different 
components. PIC methods are an important tool to understand the plasma dynamics
through solution of the electrostatic and electro-magnetic equations, see the classical introduction in \cite{birdsall1985} and  \cite{hockney1985}.
While PIC is extremely application driven, the standard schemes as discussed in
\cite{hockney1985} has to be adapted for the different problem, e.g.
strong external magnetic field, see \cite{pata2009} and \cite{speiter1999}, electro-static or electro-magnetic applications, see \cite{mark2010} and \cite{chen2011}.

Further the different numerical ideas to solve the time-dependencies,
e.g. explicit or implicit schemes, has an extremely influence to the 
numerical stability of the PIC codes.
For example, an explicit solver with integration time step $\Delta t$ 
had to satisfy $\omega_{Pe} \Delta t < 2$, where $\omega_{Pe}$ is the electron 
plasma frequency, see \cite{hockney1985} or the Courant-Friedlich-Levy condition
$c \le \frac{\Delta x}{\Delta t}$, where $c$ is the wave phase velocity.
Further the grid solver has to satisfy the electron Debye length $\Delta x \le \xi \lambda_{De}$, where $\lambda_{De}$ is the Debye length and $\xi$ is a constant of order $1$.

In the following three parts of the PIC scheme are involved and given as:
\begin{itemize}
\item Pusher (scheme to solve the mesh-free equation of motions). 
\item Solver (scheme to solve the mesh-based potential equations).
\item Interpolation (Approximation schemes to couple the mesh-free parameters with the mesh parameters)
\end{itemize}
All three parts are important and we have to deal with their numerical
error to reduce the full error of the PIC approximation.
Second, the physical constraints, as conservation of mass, momentum and
energy are important to the physical experiments and 
should be conserved by the underlying schemes.

We concentrate on improving the pusher part (time-integrator) 
to reduce the errors of the particle tracking and to optimize the 
computational amount.

The present paper is organized as follows.
In Section \ref{model}, we discuss the mathematical model.
The different time-integrator methods are presented in 
Section \ref{integrator}.
In Section \ref{num}, we present the numerical approaches of the 
electro-static model and in Section \ref{num2} we discuss the
benefits of improved integrators to the electro-magnetic model.
The conclusions are given in Section \ref{concl}.

\section{Mathematical Model} 
\label{model}

In the following, we present the two models, which are
numerically analyzed in the experiments with different 
time-integrator methods.

\subsection{Electro-Static Model} 

For the electro-static model, we deal with the
following equations.

The equations of motions (microscopic scales) are given for the
electro-static Model as:
\begin{eqnarray}  \label{poisson}
  \vec{x}' &=& \frac{d \vec{x}}{dt} = \vec{v}, \\
  \vec{v}' &=& \frac{q}{m} \vec{E} ,
\end{eqnarray}
where $\vec{x}, \vec{v}$ is the position and velocity of a particle,
$q$ is the charge and $m$ is the mass of the particle.

Further we have the additional electro-static equations 
which are solved on the grid (macroscopic scales)
\begin{eqnarray} 
\label{ele_1}
&& \nabla \cdot \vec{E} = \frac{\rho}{\epsilon_0} , \\
&&  \vec{E} = - \nabla \phi ,
\end{eqnarray} 
where $E$ is the electric potential and $E$ is the electric field, which are
given self-consistent, see \cite{hockney1985}.

While the microscopic scales are solved with ODE solvers, 
and spatial mesh-free methods. 

The idea of the electrostatic PIC is instead of calculating the 
equation of motion for each charged particle, one can solve it for so 
called "super particles", representing thousands of real particles. 

Further the solvers for the Poisson equations are based on a grid.
For such equations we have also very fast grid based solvers, e.g. iterative schemes as ILU or SuperLU.

To connect the particles, which are moving grid-less to a grid and vice versa.  We have to define an interpolation, between the particles and the grids.

We have to define a interpolation function $S(x_i - X_j)$ where $X_j$ is the grid-point and $x_i$ the position of the super-particle $i$.

With this interpolation function the weighting of the particle charge $q_i$ of the particle $i$ at the the position $x_i$ to the  grid-point $j$ gives the particle charge density at this grid point $\rho_j$
\begin{eqnarray} \label{poisson}
 \rho_j &=& \Delta x ^{-1} \sum_{i=1}^N q_i S(x_i -X_j).
\end{eqnarray}
With the cell size $\Delta x$ and $N$ the number of particles.
Further the time-steps have to resolve the Langmuir wave propagation
and the cell size have to resolve the electron Debye-length.

\subsection{Electro-Magnetic Model}

For the electro-magnetic model, we deal with the following
equations.
The equations of motions (microscopic scale) are given as:
\begin{eqnarray}  
\label{electro_mag_1}
  \vec{x}' &=& \frac{d \vec{x}}{dt} = \vec{v}, \\
\label{electro_mag_2}
  \vec{v}' &=& \frac{q}{m} (\vec{E} + \vec{v} \times \vec{B}) ,
\end{eqnarray}
where we have additional the magnetic field $\vec{B}$.

Further the electro-magnetic equations, which are solved in the PIC scheme 
on the grid (macroscopic scales) are given as:
\begin{eqnarray} 
\label{mag_1}
&&  \nabla \cdot B = 0, \\
&& \nabla \cdot E = \frac{\rho}{\epsilon_0} , \\
&& \frac{\partial E}{\partial t} = \frac{1}{\mu_0 \epsilon_0} \nabla \times B - \frac{1}{\epsilon_0} j \\
\label{mag_2}
&& \frac{\partial B}{\partial t} = - \nabla \times E ,
\end{eqnarray}

For the electromagnetic PIC (EMPIC), we have to extend
the classical electrostatic PIC, with respect to the 
magnetic fields. The particles in the field have additional
influence due to the magnetic fields, means the trajectories
are also influences with the magnetic field.

Here, we have to taken into account an extension of the 
standard integrators, see \cite{pata2009}.

\section{Time Integrators}
\label{integrator}

In the following, we discuss the different ideas:
\begin{itemize}
\item Explicit schemes,
\item Implicit schemes,
\item Hamiltonian-based schemes.
\end{itemize}

\subsection{Explicit Schemes}

Explicit schemes are forward schemes, while need not additional 
steps, e.g. inversion, to obtain the next solution, see \cite{geiser_2013}.

By the way, due to the simple algorithms, the explicit integrators have 
the drawback of the restrictions to the time-steps.
Especially for the PIC method, we have taken into account:
\begin{itemize}
\item Debye-length (smallest spatial):
$\Delta x \le  \lambda_{De}$ , $\lambda_{De}$ Debye-length. 
\item Langmuir frequency (smallest time-steps): $\Delta t \le  \frac{1}{\omega_{p}}$, $\omega_{p}$ plasma frequency.
\end{itemize}
Such restrictions neglect large time-steps and reduce the effectively
in the numerical computations.  
In the following Example \ref{exam_1}, we present a simple explicit integrator
applied to an electro-magnetic model.

\begin{example}
\label{exam_1}
Explicit Boris integrator, see \cite{chen2011}:
\begin{eqnarray}  \label{impl_1}
  \vec{x}_p^{n+1} &=&  \vec{x}_p^{n} +  \vec{v}_p^{n+1/2} \Delta t , \\
  \vec{v}_p^{n+1} &=&  \vec{v}_p^{n} +  \frac{q_s  \; \Delta t}{m_s} \left( \vec{E}_p^{n + \theta}(\vec{x}_p^{n+1/2}) + \vec{v}_p^{n+1/2} \times \vec{B}_p^{n}(\vec{x}_p^{n+1/2}) \right) ,
\end{eqnarray}
where the intermediate solutions $\vec{x}_p^{n+1/2}$ and $\vec{v}_p^{n+1/2}$ are computed a separated third step (see Strang-splitting in Equation (\ref{stran_1_1})-(\ref{stran_1_2}).
\end{example}

\subsection{Implicit Time integrators}

The idea of the implicit integrators,
are to remove the need to resolve such 
small scales (Debye length and Langmuir frequency).

The unsolved scales are kept in an approximate way allowing the 
coupling and energy transfer with larger and slower scales
that are instead fully resolved.
\begin{itemize}
\item Direct implicit methods, see \cite{welch2004}.
\item Implicit moment methods, see \cite{noguchi2007}.
\end{itemize}

The implicit schemes have the advantage of larger time-steps.
They are important in the electro-magnetic models,
while the light-wave propagation is important.
For the explicit schemes we have the restriction 
$c \Delta t \le \Delta x$ of the stability, while for the implicit schemes
it is not necessary, see \cite{welch2004}. 

\begin{example}
We deal with the following Boris-integrator,
which is done implicit, see \cite{vu1995}:
\begin{eqnarray}  \label{impl_1}
  \vec{x}_p^{n+1} &=&  \vec{x}_p^{n} +  \frac{\vec{v}_p^{n+1} + \vec{v}_p^{n} }{2} \Delta t , \\
  \vec{v}_p^{n+1} &=&  \vec{v}_p^{n} +  \frac{q_s  \; \Delta t}{m_s} \left( \vec{E}_p^{n + \theta}(\vec{x}_p^{n+1/2}) + \frac{\vec{v}_p^{n+1} + \vec{v}_p^{n} }{2} \times \vec{B}_p^{n}(\vec{x}_p^{n+1/2}) \right)  \\
    && -  \mu \Delta t  \nabla ||B_p^{n}(\vec{x}_p^{n+1/2})||  , \nonumber \\
&&  \mu = \frac{(\vec{v}_p^{n+1} - \vec{v}_p^{n}) - (\vec{v}_p^{n+1} - \vec{v}_p^{n})\vec{B}_p^{n}(\vec{x}_p^{n+1/2}) \vec{B}_p^{n}(\vec{x}_p^{n+1/2})  }{8 ||B_p^{n}(\vec{x}_p^{n+1/2})||} .
\end{eqnarray}
\end{example}

\begin{remark}
To conserve the physical behavior also with implicit methods,
we have to be taken into account, that in the larger time-step, the underlying
fluctuations are frozen or at least static.
For example, we neglect the small scales,  example quantum mechanical effects,
e.g. de Haas-van Alphen effect.
de Haas–van Alphen effect, which is a quantum mechanical effect in which the magnetic moment of a pure metal crystal oscillates as the intensity of an applied magnetic field B is increased.
The period, when plotted against $1/B$, is inversely proportional to the area 
$S$ of the extremal orbit of the Fermi surface, in the direction of the applied field.
\begin{eqnarray}
    \Delta \left( \frac{1}{B} \right) = \frac{2 \pi e}{\hbar c S}
\end{eqnarray}
where $S$ is the area of the Fermi surface normal to the direction of $B$.
\end{remark}

\subsection{Hamiltonian-based Methods}

The Hamiltonian-based methods are taken into account the conservation 
of the symplecticity with respect to the Hamiltonian-form.
Such orbit integrator methods are based on the ideas to reformulate the 
equation of motion in electric and magnetic fields, see \cite{cohen1982}.
Based on the reformulation of the equation of motions (\ref{electro_mag_1})-(\ref{electro_mag_2}) 
to a Hamiltonian form:
\begin{eqnarray}
\frac{\partial \vec{q}}{\partial t} = \frac{\partial H}{\partial \vec{p}} \\
\frac{\partial \vec{p}}{\partial t} = - \frac{\partial H}{\partial \vec{q}} \\
\end{eqnarray}
where $\vec{q} = \vec{x}$ and $\vec{p} = \vec{v}$.
Further in the case of a charged particle in an electromagnetic field,
we have also to derive the Hamiltonian.

\section{Numerical Experiments for the Electro-Static Models}
\label{num}

In the following, we apply the different time-integrator methods
with respect to their efficiency and their accuracy of the numerical 
results. We test:
\begin{itemize}
\item Euler-Forward Integrators (A-B Splitting),
\item Boris Integrator (Strang-Splitting),
\item Boris Integrator with Space Filters (Energy conserved method).
\end{itemize}

We apply an electrostatic PIC-code, which is programmed in OCTAVE.

We apply the following particle model:

1.) The trajectories of the particles are given as (microscopic scale):
\begin{eqnarray}  \label{poisson}
  \vec{x}' &=& \frac{d \vec{x}}{dt} = \vec{v}, \\
  \vec{v}' &=& \frac{q}{m} \vec{E} .
\end{eqnarray}

2.) The electrostatic field (macroscopic scale) is given as
\begin{eqnarray} 
&&  \nabla \cdot \nabla \phi = \frac{\rho}{\epsilon_0}, \\
&& \vec{E} = \nabla \phi , \\
&& \nabla \cdot \vec{E} = \frac{\rho}{\epsilon_0} .
\end{eqnarray} 
The micro- and macroscopic equations are coupled via the
approximation functions of higher order e.g., cloud-in-cell (CIC), see \cite{hockney1985}.

We apply following integrators:
\begin{itemize}
\item Euler-forward integrator (A-B splitting):
\begin{eqnarray}  \label{impl_1}
  \vec{x}_p^{n+1} &=&  \vec{x}_p^{n} +  \vec{v}_p^{n} \Delta t , \\
  \vec{v}_p^{n+1} &=&  \vec{v}_p^{n} +  \frac{q_s  \; \Delta t}{m_s}  \vec{E}_p(\vec{x}_p^{n+1})   
\end{eqnarray}
\item Boris integrator (Strang-splitting):
\begin{eqnarray}  \label{impl_1}
  \vec{x}_p^{n+1} &=&  \vec{x}_p^{n} +  \frac{\vec{v}_p^{n+1} + \vec{v}_p^{n} }{2} \Delta t , \\
  \vec{v}_p^{n+1} &=&  \vec{v}_p^{n} +  \frac{q_s  \; \Delta t}{m_s}  \vec{E}_p(\vec{x}_p^{n+1/2})   
\end{eqnarray}
or given as:
\begin{eqnarray}  \label{impl_1}
  \vec{x}_p^{n+1/2} &=&  \vec{x}_p^{n} +  \frac{\vec{v}_p^{n}}{2} \Delta t , \\
  \vec{v}_p^{n+1} &=&  \vec{v}_p^{n} +  \frac{q_s  \; \Delta t}{m_s}  \vec{E}_p(\vec{x}_p^{n+1/2})   \\
 \vec{x}_p^{n+1} &=&  \vec{x}_p^{n+1/2} +  \frac{\vec{v}_p^{n+1}}{2} \Delta t ,
\end{eqnarray}
\begin{eqnarray}  \label{impl_1}
   \vec{E}_p(\vec{x}_p^{n+1/2}) &=& \sum_i \frac{\vec{E}_i^n + \vec{E}_i^{n+1}}{2} S( \vec{x}_i - \vec{x}_p^{n+1/2} ) ,
\end{eqnarray}
where $S$ is a spline function, e.g. first order.

\end{itemize}

The PIC problem is given as:
$q_s = -1$, $m_s = 1.0$ with periodic boundary conditions. We apply the first order approximation to the two stream instability , see \cite{lapenta2011}.

To see the development of the integrators it is sufficient to use the max-norm with a reference result of very fine Euler-integrator result: 
\begin{eqnarray}  \label{impl_1}
 err_{method} = \max_{p=1, \ldots, P} \| x_{p,method} - x_{p, reference} \| ,
\end{eqnarray}
where $\mbox{method} = \{\mbox{Euler}, \mbox{Strang}\}$
and $\mbox{reference}$ is the Euler-method with very fine resolutions.

We start with a converged Euler solution which is given with $x_{p, reference (\Delta t_{fine})}$, while the difference of $\max_{p=1, \ldots, P} \| x_{p, reference (\Delta t_{fine})} - x_{p, reference (\Delta t_{fine}/2)} \| \le \max_{p=1, \ldots, P} \| x_{p, method(\Delta t)} - x_{p, reference (\Delta t_{fine}/2)} \|$, means the reference solution is numerically converged.

This means: we set the Euler solution as benchmark solution and compare the other integrators with it.

At first we need the benchmark result of the very fine Euler-integrator.
The result is given in the Figure \ref{euler_1}:

\begin{figure}[ht]
	\centering
	\includegraphics[width=\textwidth]{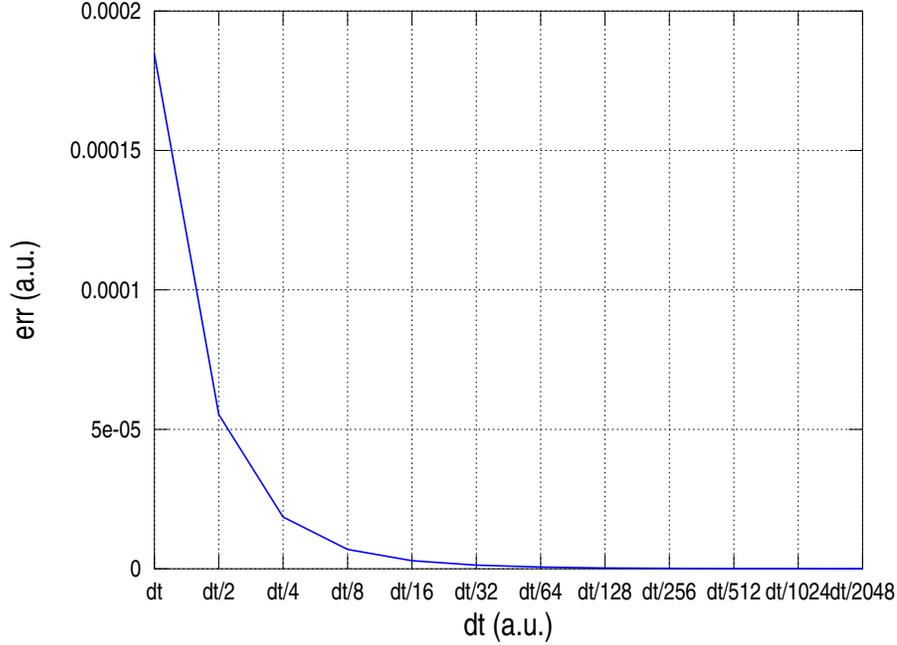}       
	\caption{\label{euler_1} Reference solution based on the Euler integrator with sufficient small time-steps, err$_{euler}$ for various $\Delta$t's.}
\end{figure}

%

Now we compare the convergence behavior of the different improved 
integrators and underlie the Euler integrator benchmark solution
with related integrator.

The improved integrators are given in the following Figure \ref{boris_1},
with the max-Norm and a fine resolved reference solution done with 
Euler method.
\begin{figure}[ht]
	\centering
                \includegraphics[width=\textwidth]{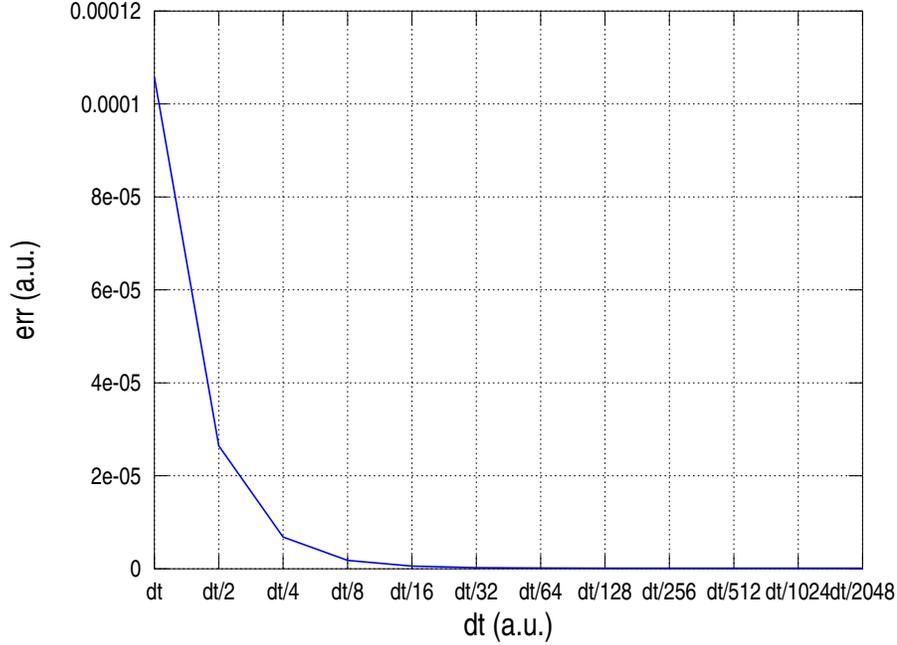}
	\caption{\label{boris_1} Boris Integrator solutions for various $\Delta$t; the reference solution is done Euler Integrator and $\Delta$t/512.}
\end{figure}

As you can see the Boris-integrator converges faster than the Euler-integrator. The question is if we can improve the convergence further due to filtering of the electric field.

In the next figure \ref{borisfilter}, we present an improvement
to the Boris-integrator based on the filter-technique.
Based on the filtering, we could improve the convergence results of the
integrators and obtain fast convergent results with time-steps about $\Delta t/8$. Due to improvement of the approximation error, we have a shift to
our standard reference solution done by Euler-integrators.
If we shift our results we obtain the improved convergence rates of such 
novel methods. 
\begin{figure}[ht]
	\centering
	        \includegraphics[width=\textwidth]{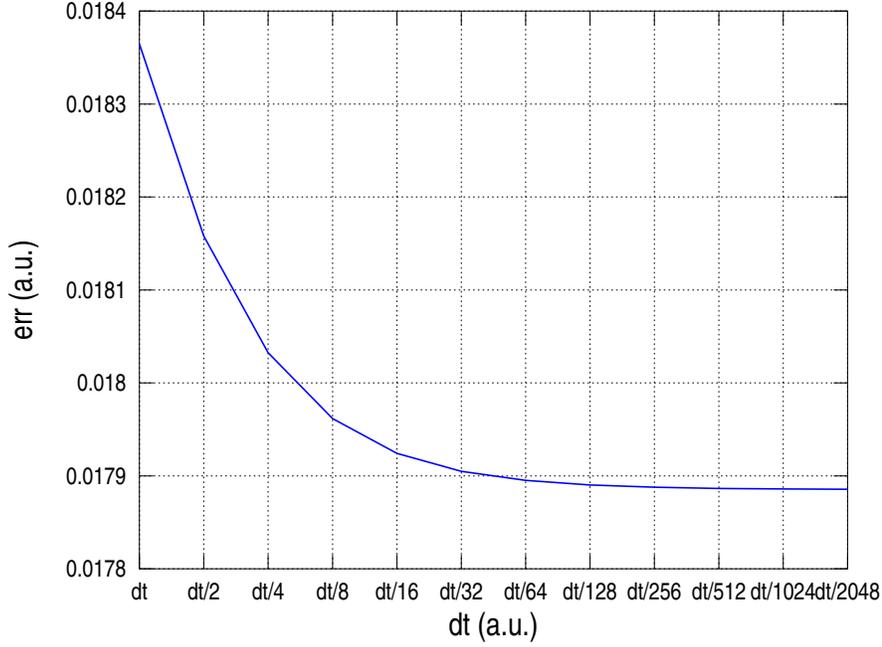}
	 \label{fig:boris_filter}       
	\caption{\label{borisfilter} The improvement based on the filter techniques for the Boris integrator, we see an offset, based on the reference solution of the Euler integrator.}
\end{figure}

\section{Numerical Experiments for the Electro-Magnetic Models}
\label{num2}

In the following, we apply the different time-integrator methods
with respect to their efficiency and their accuracy of the numerical 
results. We test:
\begin{itemize}
\item Euler-Forward Integrators (A-B Splitting),
\item Boris Integrator (Strang-Splitting),
\item Cyclotronic Integrator (Symplectic method).
\end{itemize}

We apply an electrostatic PIC-code, which is programmed in OCTAVE.

We apply the following particle model:

1.) The trajectories of the particles are given as (microscopic scale):
\begin{eqnarray}  \label{poisson}
  \vec{x}' &=& \frac{d \vec{x}}{dt} = \vec{v}, \\
  \vec{v}' &=& \frac{q}{m} (\vec{E} + v \times \vec{B} ) ,
\end{eqnarray}
where $\vec{B}$ is an external constant magnetic field,
given as:
\begin{eqnarray}  \label{poisson}
  \omega_B &=& \frac{q |\vec{B}|}{m} , 
\end{eqnarray}
and we have two equations:
\begin{eqnarray}  \label{poisson}
  x_x' &=&  v_x, \\
  x_y' &=&  v_y, \\
  v_x' &=& \frac{q}{m} E_x + \omega_B v_y  , \\
  v_y' &=& \frac{q}{m} E_y - \omega_B v_x  .
\end{eqnarray}
2.) The electro field (macroscopic scale) is given as:
\begin{eqnarray} 
&&  \nabla \cdot \nabla \phi = \frac{\rho}{\epsilon_0}, \\
&& \vec{E} = \nabla \phi , \\
&& \nabla \cdot \vec{E} = \frac{\rho}{\epsilon_0} ,
\end{eqnarray} 
while here we deal with two dimensional spatial operators,
e.g. $\Delta x$ and $\Delta y$ are the spatial steps
and we have two-dimensional differential operators, see the dicretization
scheme for the Maxwell equation in Appendix \ref{appendix}.

The PIC-problem is given

The micro- and macroscopic equations are coupled via the
approximation functions.

The integrators are given as:
\begin{itemize}
\item Euler-forward integrator (A-B splitting):
\begin{eqnarray}  \label{impl_1}
  \vec{x}_p^{n+1} &=&  \vec{x}_p^{n} +  \vec{v}_p^{n} \Delta t , \\
  \vec{v}_p^{n+1} &=&  \vec{v}_p^{n} +  \frac{q_s  \; \Delta t}{m_s} ( \vec{E}_p(\vec{x}_p^{n+1})  +\vec{v}_p^{n} \times  \vec{B} )
\end{eqnarray}
where $\vec{v}_p^{n} \times  \vec{B} = \left( \begin{array}{c} \omega_b v_y^n \\ - \omega_b v_x^n   \end{array} \right)$

\item Boris integrator (Strang-splitting):
\begin{eqnarray}  \label{impl_1}
\label{stran_1_1}
  \vec{x}_p^{n+1/2} &=&  \vec{x}_p^{n} +  \frac{\vec{v}_p^{n}}{2} \Delta t , \\
  \vec{v}_p^{n+1/2} &=&  \vec{v}_p^{n} +  \frac{q_s  \; \Delta t}{2 m_s}  \vec{E}_p(\vec{x}_p^{n+1/2})   \\
  \vec{v}_p^{n*} &=&  \vec{v}_p^{n+1/2} +  \frac{q_s  \; \Delta t}{m_s} ( \vec{v}_p^{n+1/2} \times  \vec{B} )   \\
 \vec{x}_p^{n*+1/2} &=&  \vec{x}_p^{n} +  \frac{\vec{v}_p^{n*}}{2} \Delta t , \\
  \vec{v}_p^{n+1} &=&  \vec{v}_p^{n*} +  \frac{q_s  \; \Delta t}{2 m_s}  \vec{E}_p(\vec{x}_p^{n*+1/2})   \\
\label{stran_1_2}
 \vec{x}_p^{n+1} &=&  \vec{x}_p^{n+1/2} +  \frac{\vec{v}_p^{n+1}}{2} \Delta t ,
\end{eqnarray}
\begin{eqnarray}  \label{impl_1}
   \vec{E}_p(\vec{x}_p^{n+1/2}) &=& \sum_i \frac{\vec{E}_i^n + \vec{E}_i^{n+1}}{2} S( \vec{x}_i - \vec{x}_p^{n+1/2} ) , \\
  \vec{E}_p(\vec{x}_p^{n*+1/2}) &=& \sum_i \frac{\vec{E}_i^n + \vec{E}_i^{n+1}}{2} S( \vec{x}_i - \vec{x}_p^{n*+1/2} ) ,
\end{eqnarray}
where $S$ is a spline function, e.g. first order.
\item Cyclotronic Integrator (cyclic Splitting):
We deal with a $\vec{B} = B \vec{e}_z$ and Larmor angular frequency $\Omega= q |\vec{B}|/m$, \\

Step 1:
\begin{eqnarray}  \label{impl_1}
  z^{n+1/2} &=&  z^{n} +  v_{z}^{n} \frac{\Delta t}{2} , \\
  x^{n+1/2} &=&  x^{n} +  \frac{v_y^{n} - v_y^{n} \cos(\Omega \frac{\Delta t}{2}) + v_x^{n} \sin(\Omega \frac{\Delta t}{2})}{\Omega} , \\
  y^{n+1/2} &=&  y^{n} +  \frac{- v_x^{n} + v_x^{n} \cos(\Omega \frac{\Delta t}{2}) + v_y^{n} \sin(\Omega \frac{\Delta t}{2})}{\Omega} , \\
  v_x^{n+1/2} &=&  v_x^{n}  \cos(\Omega \frac{\Delta t}{2}) + v_y^{n} \sin(\Omega \frac{\Delta t}{2}) , \\
  v_y^{n+1/2} &=&  v_y^{n}  \cos(\Omega \frac{\Delta t}{2}) - v_x^{n} \sin(\Omega \frac{\Delta t}{2}) ,
\end{eqnarray}
Step 2:
\begin{eqnarray}  \label{impl_1}
  \vec{v}^{n*} &=&  \vec{v}^{n+1/2} + \frac{q_s  \; \Delta t}{2 m_s}  \vec{E}_p(\vec{x}^{n+1/2}), \\
 \vec{E} & = & - \nabla \phi ,
\end{eqnarray}
Step 3:
\begin{eqnarray}  \label{impl_1}
  z^{n+1} &=&  z^{n*} +  v_{z}^{n*} \frac{\Delta t}{2} , \\
  x^{n+1} &=&  x^{n*} +  \frac{v_y^{n*} - v_y^{n*} \cos(\Omega \frac{\Delta t}{2}) + v_x^{n*} \sin(\Omega \frac{\Delta t}{2})}{\Omega} , \\
  y^{n+1} &=&  y^{n*} +  \frac{- v_x^{n*} + v_x^{n*} \cos(\Omega \frac{\Delta t}{2}) + v_y^{n*} \sin(\Omega \frac{\Delta t}{2})}{\Omega} , \\
  v_x^{n+1} &=&  v_x^{n*}  \cos(\Omega \frac{\Delta t}{2}) + v_y^{n*} \sin(\Omega \frac{\Delta t}{2}) , \\
  v_y^{n+1} &=&  v_y^{n*}  \cos(\Omega \frac{\Delta t}{2}) - v_x^{n*} \sin(\Omega \frac{\Delta t}{2}) ,
\end{eqnarray}
\end{itemize}

The particle in cell problem is given with the following parameters:
$q_s = -1$, $m_s = 1.0$ and $B=1$ which gives a Larmor-frequency of 1. The boundary conditions are periodic.

The initial conditions are equal for all $x$ and $y$. The particles
are uniformly distributed in the space and the initial velocities
of the particles are given as in the 1D case for $x$ and $y$.
The initial conditions are generated via $rand()$ and are stored,
such that we have in each experiment the same initial conditions.

The errors of the pusher is given as:
\begin{eqnarray}  \label{impl_1}
 err_{L_2, method} = ( \sum_{p=1}^p \Delta x \; (x_{p,method} - x_{p, reference} )^2 )^{1/2} ,
\end{eqnarray}
where $\mbox{method} = \{\mbox{Euler}, \mbox{Strang}\}$
and $\mbox{reference}$ is the Euler-method with very fine resolutions,
the spatial step of the grid is given as $\Delta x$.

As in the 1D case start with a converged Euler solution which is given with $x_{p, reference (\Delta t_{fine})}$, while the difference of $\max_{p=1, \ldots, P} \| x_{p, reference (\Delta t_{fine})} - x_{p, reference (\Delta t_{fine}/2)} \| \le \max_{p=1, \ldots, P} \| x_{p, method(\Delta t)} - x_{p, reference (\Delta t_{fine}/2)} \|$, means the reference solution is numerically converged.

The numerical convergence is given as:
\begin{eqnarray}  \label{impl_1}
 err_{method, \Delta_{method} t} = \max_{p=1, \ldots, P} \| x_{p,method, \Delta_{method} t} - x_{p, reference (\Delta t_{fine})} \| ,
\end{eqnarray}
where $\Delta t_{fine} = \Delta t / 2048$ and $\Delta t_{method} = \{ \Delta t , \Delta t/2 , \ldots,  \Delta t/2048 \}$. 
%

The errors of the pusher is given as:
\begin{eqnarray}  \label{impl_1}
 err_{method, \Delta x, \Delta t} = \max_{p=1, \ldots, P} \| x_{p, \Delta x, \Delta t, method} - x_{p, \Delta x_{fine}, \Delta t_{fine}, reference} \| ,
\end{eqnarray}
where we assume $\Delta x = \Delta y$ and the convergence tableau is given with
the spatial scales
$\Delta x, \Delta x/2, \ldots, \Delta x/16 = \Delta x_{fine}$, $\Delta y, \Delta y/2, \ldots, \Delta y/16 = \Delta y_{fine}$ 
and the time scales $\Delta t , \Delta t/2 , \ldots,  \Delta t/2048 = \Delta t_{fine}$.

The numerical convergence rate is given as:
\begin{eqnarray}  \label{impl_1}
 \rho_{err_{method, \Delta_{method} t}} = \frac{\log(\frac{err_{method, \Delta_{method} t/2}}{err_{method, \Delta_{method} t}})}{\log(0.5)} ,
\end{eqnarray}

\begin{remark}
The convergence-rates are optimal for sufficient small time and spatial steps,
where we obtain optimal results for $\Delta x/512$ and $\Delta t/128$.
\end{remark}

Here, we applied a spatial-temporal convergence with the different
refined time and spatial steps, e.g. $\Delta x, \Delta x/2, \ldots, \Delta x/2048$, $\Delta y, \Delta y/2, \ldots, \Delta y/2048$ and  $\Delta t, \Delta t/2, \ldots, \Delta t/512$.

\begin{figure}[ht]
	\centering
                \includegraphics[width=\textwidth]{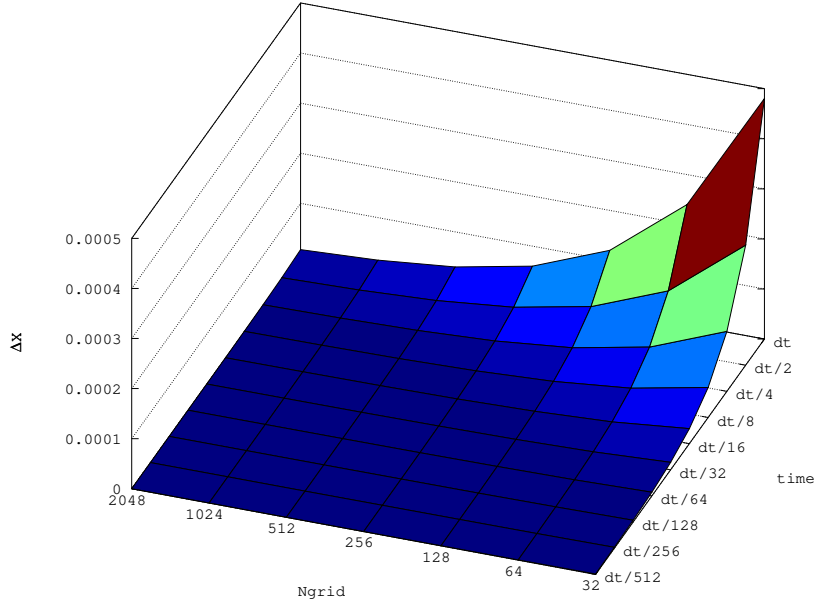}
                \label{fig:conveuler}
\caption{\label{3d_conv} Euler-Scheme: 3d space-time-convergence diagram with x-axis: $\Delta x, \Delta x/2, \ldots, \Delta x/16$, $\Delta y, \Delta y/2, \ldots, \Delta y/16$ and y-axis: $\Delta t, \Delta t/2, \ldots, \Delta t/16$, and z-axis is the error.}
\end{figure}

The results are given in the Figure \ref{3d_conv}.

The errors of the pusher is given as:
\begin{eqnarray}  \label{impl_1}
 err_{method, \Delta x, \Delta t} = \max_{p=1, \ldots, P} \| x_{p, \Delta x, \Delta t, method} - x_{p, \Delta x_{fine}, \Delta t_{fine}, reference} \| , \\
 err_{L_2, method, \Delta x, \Delta t} = ( \sum_{p=1}^p \Delta x \; (x_{p, \Delta x, \Delta t, method} - x_{p, \Delta x_{fine}, \Delta t_{fine}, reference} )^2 )^{1/2} ,
\end{eqnarray}
where we assume $\Delta x = \Delta y$ and the convergence tableau is given with
the spatial scales
$\Delta x, \Delta x/2, \ldots, \Delta x/16 = \Delta x_{fine}$, $\Delta y, \Delta y/2, \ldots, \Delta y/16 = \Delta y_{fine}$ 
and the time scales $\Delta t , \Delta t/2 , \ldots,  \Delta t/1024 = \Delta t_{fine}$. 

\begin{remark}
The convergence-rates are optimal for sufficient small time and spatial steps,
where we obtain optimal results for $\Delta x/512$ and $\Delta t/128$.
\end{remark}

\begin{figure}[ht]
	\centering
             \includegraphics[width=\textwidth]{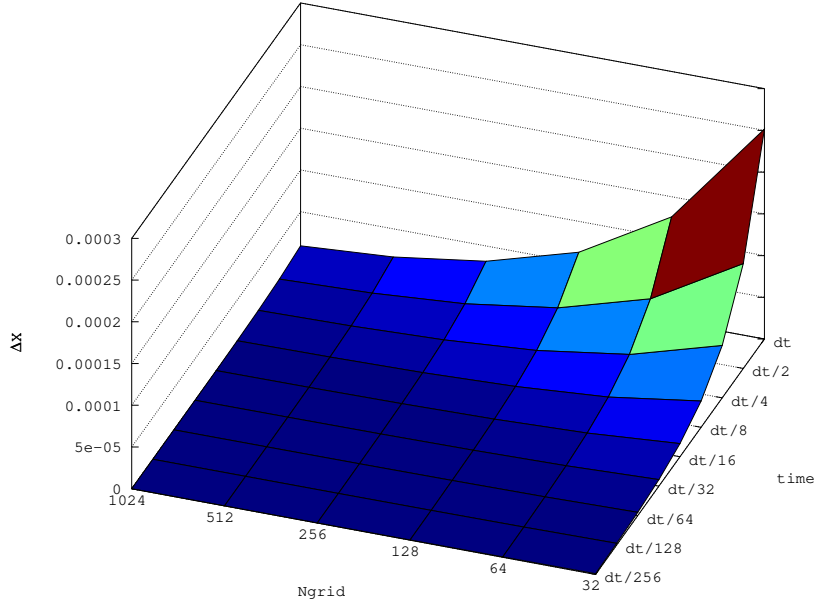}
\caption{\label{boris_2d} Boris-Integrator: 3d space-time-convergence diagram with x-axis: $\Delta x, \Delta x/2, \ldots, \Delta x/16$, $\Delta y, \Delta y/2, \ldots, \Delta y/16$ and y-axis: $\Delta t, \Delta t/2, \ldots, \Delta t/16$, and z-axis is the error.}
\end{figure}

A next improvement of the integrators is given with the 
Boris-integrator, which is a second order scheme.
The results of the Boris-integrator is given in the Figure \ref{boris_2d}.

\begin{figure}[ht]
	\centering
             \includegraphics[width=\textwidth]{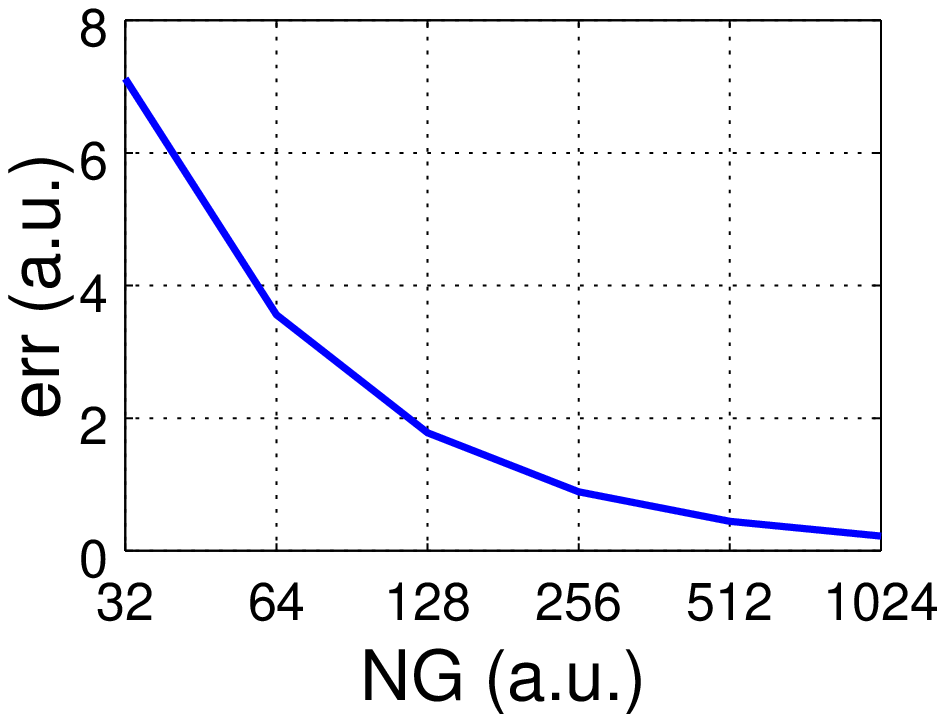}
\caption{\label{boris_2d_ref} Boris-Integrator: convergence diagram with  a finer resolution of the errors: x-axis: $\Delta x, \Delta x/2, \ldots, \Delta x/16$, $\Delta y, \Delta y/2, \ldots, \Delta y/16$ and y-axis s the error.}
\end{figure}

The results of the Boris-integrator minus the reference is given in the Figure \ref{boris_2d_ref}. As a reference the result of the Euler integrator with $NG=2048$ and $DT=dt/1024$ is used.
Because the improvement over DT is very small compared to the improvement over NG the averages over DT results were taken and plotted against NG.

\begin{figure}[ht]
	\centering
             \includegraphics[width=\textwidth]{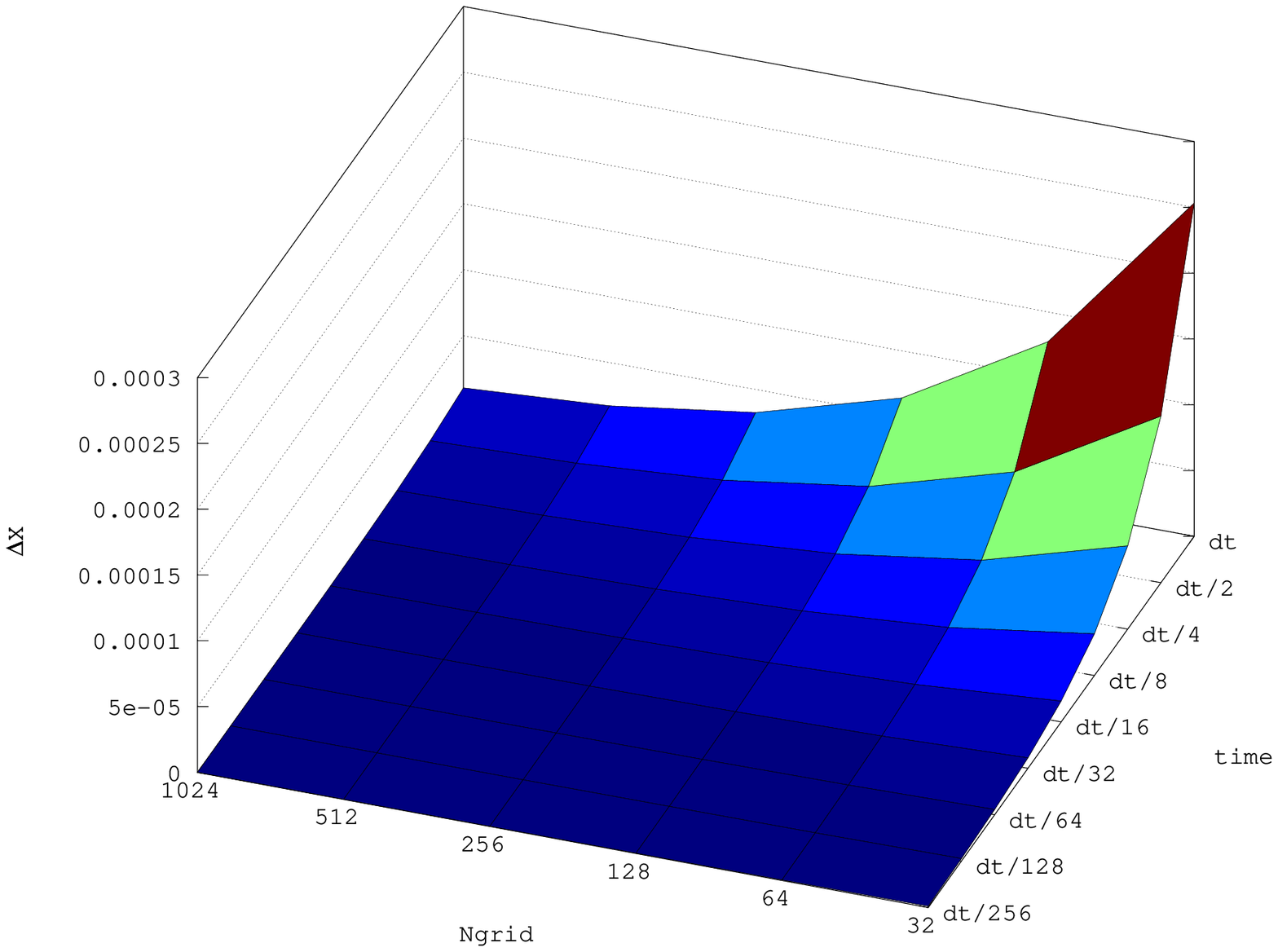}
\caption{\label{boris-filter_2d} Boris-Filter-Integrator: 3d space-time-convergence diagram with x-axis: $\Delta x, \Delta x/2, \ldots, \Delta x/16$, $\Delta y, \Delta y/2, \ldots, \Delta y/16$ and y-axis: $\Delta t, \Delta t/2, \ldots, \Delta t/16$, and z-axis is the error.}
\end{figure}

As shown in the 1D electrostatic case the Boris-Filter-integrator should give further improvement.
The results of the Boris-Filter-integrator is given in the Figure \ref{boris-filter_2d}.

\begin{figure}[ht]
	\centering
             \includegraphics[width=\textwidth]{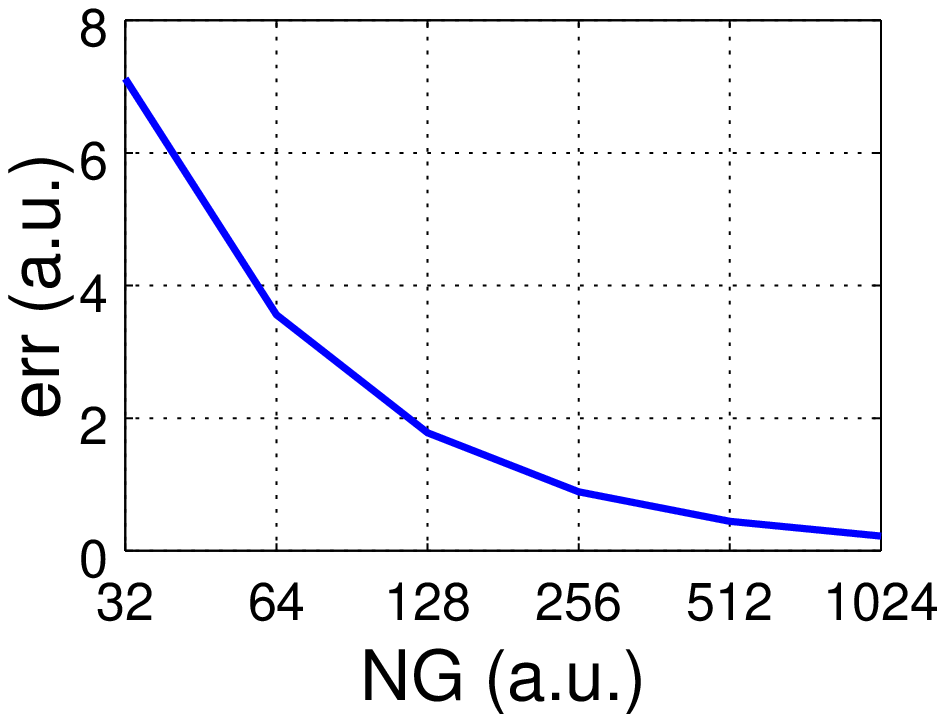}
\caption{\label{boris-filter_2d_ref} Boris-Filter-Integrator: convergence diagram with  a finer resolution of the errors: x-axis: $\Delta x, \Delta x/2, \ldots, \Delta x/16$, $\Delta y, \Delta y/2, \ldots, \Delta y/16$ and y-axis is the error.}
\end{figure}

The results of the Boris-Filter-integrator minus the reference (same as for Boris-Integrator) is given in the Figure \ref{boris-filter_2d_ref}.

This time the improvement of the Boris-filter-integrator is much smaller than in 1D case. This could be because both Boris-integrators were tested against time and not against $\Delta x$. As a result the Boris-filter-integrator only gives improvement in terms of time-steps.
The next integrator is not hardly used in plasma simulation but should give much more improvements than the Boris-filter-integrator.

The cyclotronic-integrator handles the magnetic-field different because it is Hamilton based. The Larmor-frequency goes into the rotation-operator and should be much more precise in terms of handling the magnetic-field.
The results of the Cyclotronic integrator is given in the Figure \ref{cyclo_2d_ref}.
\begin{figure}[ht]
	\centering
             \includegraphics[width=\textwidth]{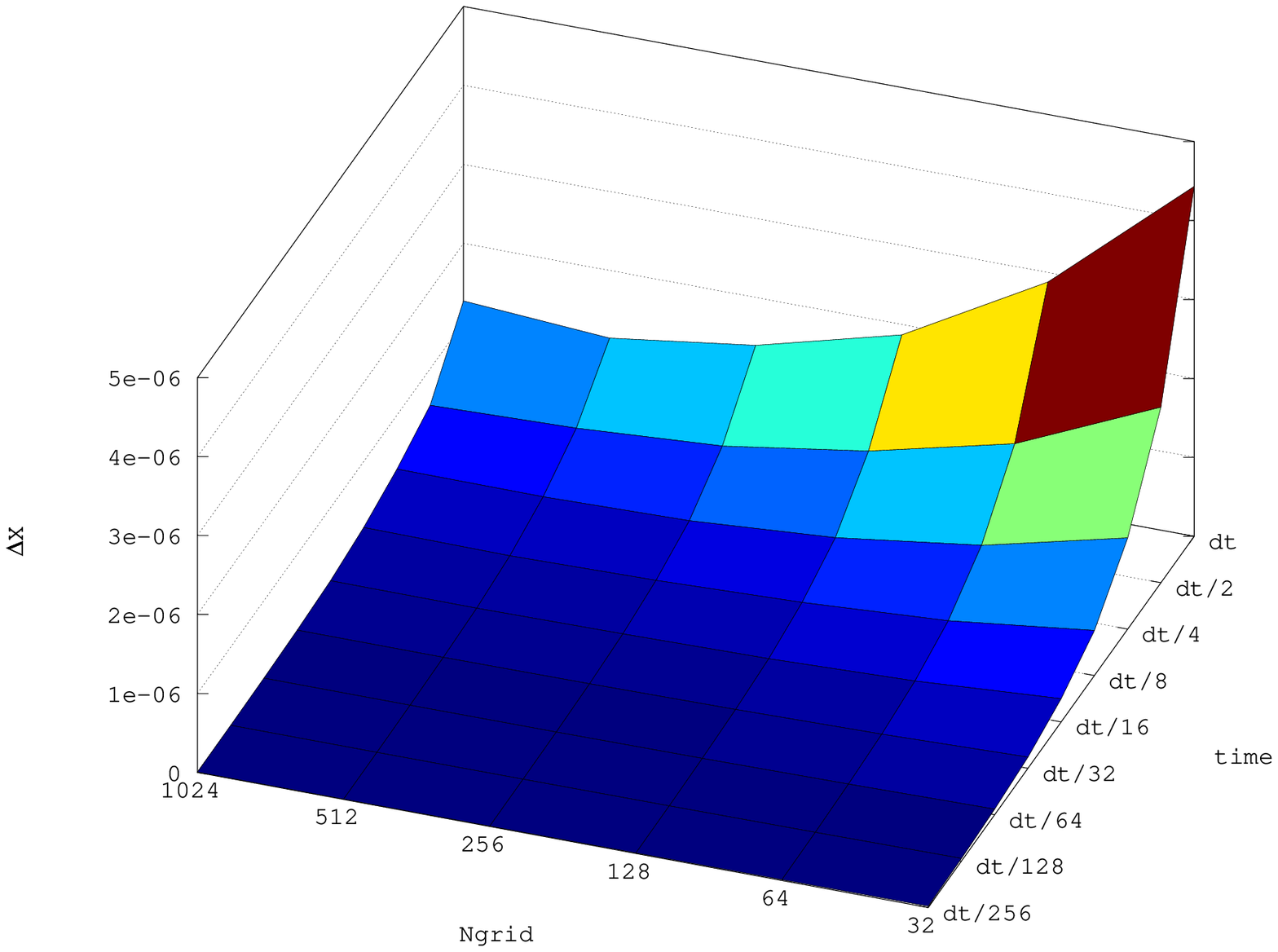}
\caption{\label{cyclo_2d_ref} Cyclotronic-Integrator: 3d space-time-convergence diagram with  a finer resolution of the errors: x-axis: $\Delta x, \Delta x/2, \ldots, \Delta x/16$, $\Delta y, \Delta y/2, \ldots, \Delta y/16$ and y-axis: $\Delta t, \Delta t/2, \ldots, \Delta t/16$, and z-axis is the error.}
\end{figure}

It can be seen that the error is much lower from the coarsest time- and grid-resolutions. This is a further improvement.

Figure \ref{cyclo_2d_ref} shows the results of the Cyclotronic integrator minus the reference (same as for Boris-Integrator). It can be seen that the difference between the cyclic-integrator and the converged Euler-solution is greatly improved.
\begin{figure}[ht]
	\centering
             \includegraphics[width=\textwidth]{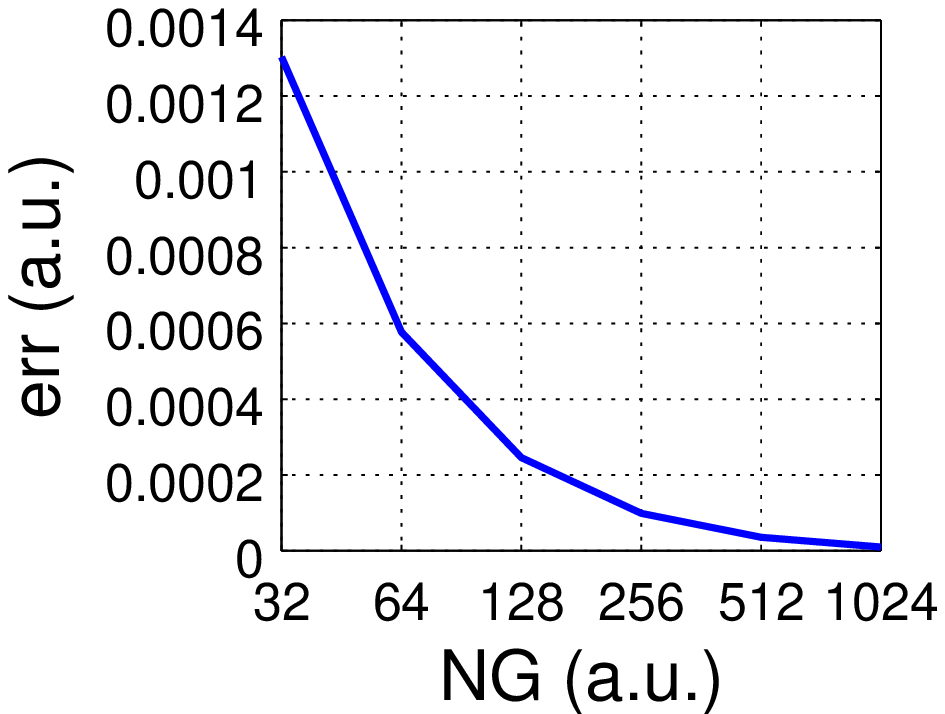}
\caption{\label{cyclo_2d_ref} Cyclotronic-Integrator: convergence diagram with  a finer resolution of the errors: x-axis: $\Delta x, \Delta x/2, \ldots, \Delta x/16$, $\Delta y, \Delta y/2, \ldots, \Delta y/16$ and y-axis is the error.}
\end{figure}

\begin{remark}
In the convergence rates, we see the benefits of the orbit integrator methods,
they decrease the error optimal, while they conserve the symplecticity of the
method. Based on the application of the filters, the conservation of the energy is obtained but the numerical error is not decreased as in the cyclotronic integrator.
Here the symplecticity is sufficient for the momentum conservation and
the long time energy conservation.
\end{remark}

 \begin{remark}

The run-time analysis is performed on an Intel Core2Quad CPU Q9400 @ 2.66GHz × 1 with 4GB of RAM. The operating system was Ubuntu 12.04LTS (64bit, Kernel:3.8.0-38-generic) and octave 3.8.1 was used to run the pic-codes. The time-step used was $DT=dt/1$. As a result you can see that the Euler and cyclic run-times develop similar. The same for Boris and Boris-filter. The Boris and Euler run-times do not develop parallel. The slopes for the Boris-integrators are slightly higher than for Euler and cyclic-integrators.
\begin{figure}[ht]
	\centering

             \includegraphics[width=\textwidth]{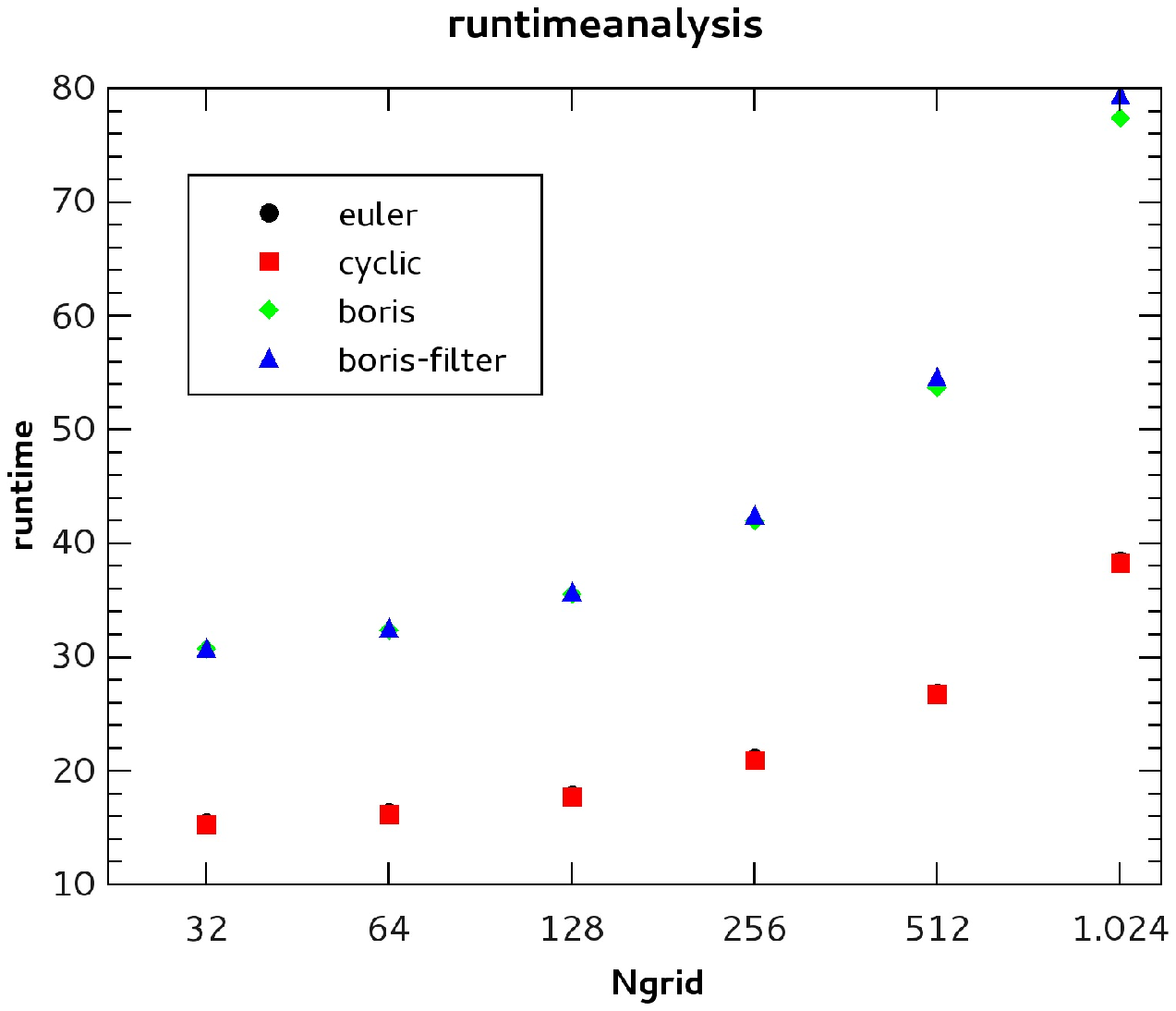}
\caption{\label{runtime} Run-time analysis of the four integrators: Euler, Boris, Boris-filter, cyclic. Time-step used: $dt/1$; Ngrid from 32 to 1024. }
\end{figure}

As a result it is recommended to use the Cyclotronic Integrator in case of electromagnetic problems because it improves the run-times of electromagnetic pic-code and converges faster.

\end{remark}

\section{Conclusion}
\label{concl}

In this paper, we discussed the benefits of the
different time-integrators for the different PIC schemes.
General integrators are at least 
Boris-integrator (explicit) or direct implicit methods (implicit).
While the Boris-integrator leaks to very small time-steps,
the implicit methods have their drawback in implicit
handling of the equations (inverse problems).
The best results are obtained by orbit integrators, while preserving the
constraints and we could apply larger time-steps.
In case of 2D electromagnetic problems the Cyclotronic Integrator has improved convergence and run-times similar to Euler. Therefor it is superior to the Boris integrator with filtering.
The electron self-force was not investigated in this paper.
In future, we discuss the numerical analysis of the improved
orbit integrators with respect to the PIC schemes.

\section{Appendix}
\label{appendix}

For the time-discretization of Maxwell-equation (\ref{mag_1})-(\ref{mag_2}), we apply
a fractional stepping scheme, see \cite{inno2013}, which is given as:
\begin{eqnarray}
\label{max}
\nabla \times E^{n + \theta} + \frac{B^{n+1} - B^n}{\Delta t} = 0 , \\
\nabla \times B^{n + \theta} - \frac{1}{c^2} \frac{E^{n+1} - E^n}{\Delta t} = \mu_0 J^{n + 1/2} , \\
\nabla \cdot E^{n + \theta} =  \frac{\rho^{n + \theta}}{\epsilon_0} , \\
\nabla \cdot E^{n + \theta} =  \nabla \cdot E^{n + \theta} ,
\end{eqnarray}
and the generic quantity $\phi$ at time $n + \theta$ is given as
$\phi^{n + \theta} = \theta \phi^{n+1} + (1 - \theta) \phi^n$.
For $\theta = 0.5$, we obtain a second order scheme (Crank-Nicolson
scheme), see \cite{hockney1985}.

%

\bibliographystyle{plain}

\end{document}